\documentclass[11pt]{amsart}
\usepackage{graphicx}
\usepackage{amssymb}
\usepackage{epstopdf}

\usepackage{amsmath,amsfonts,amsthm,amsopn,cite,mathrsfs}

 \theoremstyle{theorem}

\theoremstyle{remark}

\newcommand{\CC}{\mathbb{C}\,}

\newcommand{\RR}{\mathbb R}
\newcommand{\ZZ}{\mathbb Z}

\newcommand{\Lt}{{L^2(\RR)}}
 \newcommand{\PW}{PW_{[0, \beta]}}
  \newcommand{\PWm}{PW_{[0, \beta(M)]}}

\newcommand{\beq}{\begin{equation}}
 \newcommand{\eeq}{\end{equation}}
\newcommand{\beqq}{\begin{equation*}}
 \newcommand{\eeqq}{\end{equation*}}

\newcommand{\Ga}{{\mathcal G}(\Lambda,M) }

\begin{document}

\sloppy
\title[Irregular Gabor frames of Cauchy kernels ]
{Irregular Gabor frames of Cauchy kernels }

\author{Yurii Belov}
\address{Yurii Belov,
\newline Department of Mathematics and Computer Sciences, 
%\newline 
St.~Petersburg State University,  29 Line 14th (Vasilyevsky Island), 199178, St. Petersburg, Russia,
\newline {\tt j\_b\_juri\_belov@mail.ru} }
\author{Aleksei Kulikov}
\address{Aleksei Kulikov,  
\newline Department of Mathematics and Computer Sciences,
%\newline 
St.~Petersburg State University,  29 Line 14th (Vasilyevsky Island), 199178 St. Petersburg, Russia,
\newline Department of Mathematical Sciences, 
%\newline 
Norwegian University of Science and Technology, NO-7491 Trondheim, Norway,
\newline {\tt lyosha.kulikov@mail.ru} }
\author{Yurii Lyubarskii}
\address{Yurii Lyubarskii,  
\newline Department of Mathematics and Computer Sciences
%\newline 
St.~Petersburg State University,  29 Line 14th (Vasilyevsky Island), 199178, St. Petersburg, Russia,
\newline Department of Mathematical Sciences, 
%\newline 
Norwegian University of Science and Technology, NO-7491 Trondheim, Norway,
\newline {\tt yuralyu@gmail.com} }

%\thanks{The work was supported by the Russian Science Foundation grant 19-11-00058 and by Project 275113 of the Research Council of Norway.} 

\thanks{This work was carried out in the framework of the project 20-61-46016 by the Russian Science Foundation}

\maketitle

 In this note we discuss  the frame property in $\Lt$ of the irregular Gabor system generated by a single Cauchy kernel,
 i.e., the system
 \beqq
 \label{1}
 \Ga= \{\phi_{\lambda,\mu}(t)\}_{\lambda\in \Lambda, \mu\in M}=\left \{ \frac{e^{-2i\pi \mu t}}{t-\lambda-iw} \right \}_{\lambda\in \Lambda, \mu\in M },
 \eeqq
 here $w\in \CC$, $\Re w \neq 0$,  and $\Lambda=\{\lambda\} \subset \RR$, $M=\{\mu\} \subset \RR$ are, generally speaking,  irregular sets, their properties will be described further in this note. 
 In other words, we look for the  frame inequality
 \beqq
 \label{2}
  A\|f\|^2 \leq \sum_{\lambda\in \Lambda, \mu \in M} |\langle f, \phi_{\lambda,\mu}\rangle |^2 \leq B\|f\|^2, \ f\in \Lt.
 \eeqq
 
 The reason we wrote this note is twofold. First, in contrast to the  (now) classical  rectangular lattices $\alpha\ZZ\times\beta\ZZ$, not much is known
 about irregular ones $\Lambda \times M$. The recent breakthrough related to semiregular lattices of the form $\Lambda\times \beta \ZZ$ has been achieved  in \cite{GRS},  where the authors considered the  Gabor
 frames, generated by Gaussian totally positive functions of finite type. We also refer the reader to \cite{GRS} for the history of the problem.  Unfortunately, we are not able to apply the techniques from \cite{GRS}
 to  $\Ga$ even in the case, when $M=\beta \ZZ$ for some $\beta >0$. Second, the proof we suggest is very simple, much simpler than  the known ones for rectangular lattices, see  \cite{J,GS}.
 We expect this proof can serve as a model in more general settings.
 
\medskip

We use the following notation and definition.

{\bf Definition 1.} Given $\beta >0$, by the Paley-Wiener space $\PW$ we mean
\beqq
\label{3}
\PW=\{f:\, f(z)=\int_0^\beta e^{2i\pi\xi z} \hat{f}(\xi)d\xi, \, \hat{f}\in L^2(0, \beta)\}.
\eeqq 
This space consists of entire functions of exponential type, which belong to $\Lt$ and have the indicator diagram included in $[-2i\pi\beta,0]$.
We refer the reader to \cite{L} for the detailed description of the Paley-Wiener spaces as well as for other facts on entire functions. 

\medskip

{\bf Definition 2.} We  say that the set $\Lambda\subset \RR$ is {\em sampling}  for $\PW$ if, for some constants $0<A,B<\infty$, {\em (sampling constants)}
\beq
\label{4}
A\|f\|^2
\leq \sum_{\lambda\in \Lambda}  |f(\lambda)|^2 \leq B\|f\|^2, \, f\in \PW.
\eeq
Such sets have been completely  described in \cite{OS}. Typically, they have density bigger than $1/\beta$, yet their structure may be rather complicated, in particular 
they need not have density strictly bigger than 
$1/\beta$, also they need not contain subsets which are complete interpolating sequences for $\PW$.

\medskip

{\bf Definition 3.} The set $M=\{\mu_n\}\subset \RR$, $\mu_n<\mu_{n+1}$ is {\em locally finite} if
\beq
\label{4a}
\beta(M):=\sup \{\beta_n\}< \infty, \ \text{here} \ \beta_n=\mu_{n+1}-\mu_n, 
\eeq
and
\beq
\label{4b}
\sup_{x\in \RR} \# \{M\cap [x,x+1]\} < \infty.
\eeq

\medskip

%{\bf Theorem}
%{\em Let $M\subset \RR$ be a locally finite set and $\Lambda$ be a sampling set for $\PWm$. Then the system $\Ga$ is a frame in $\Lt$.}
 
 {\bf Theorem} {\em  Let $M\subset \RR$ and $  \Lambda\subset \RR$. The following statements are equivalent:
 
 1) $M$ is a locally finite set and $\Lambda$ is a sampling set for $\PWm$;
 
 2)  The system $\Ga$ is a frame in $\Lt$.}
 
 \medskip 
 
 {\em Proof  1) $\Rightarrow$ 2).} Let for definiteness $\Re w >0$ so $iw$ lies in the upper half-plane. In what follows it is convenient to
 simplify notation, by writing  $\phi_{\lambda,n}$ instead of   $\phi_{\lambda,\mu_n}$. Given $f\in \Lt$, denote
 \beq
 \label{5}
 c_{\lambda,n}=\langle \phi_{\lambda,n}, \bar{f}\rangle = \int_{-\infty}^\infty f(t)
       \frac{e^{-2i\pi\mu_n t}}{t-\lambda-iw}dt.
 \eeq
 Let
 \beqq
 \label{6}
 f_k(t)=\int_{\mu_k}^{\mu_{k+1}}\hat{f}(\xi)e^{2i\pi\xi t}d\xi,
 \eeqq
 as always, $\hat{f}$ stays for the Fourier  transform of $f$. 
 
 We have 
 \beq
 \label{7}
 f(t)=\sum_k f_k(t), \ \|f\|^2=\sum_k \|f_k\|^2.
 \eeq
  
Denote  also 
\beqq
 h_k(t):= e^{-2i\pi \mu_kt}f_k(t)\in \PWm. 
\eeqq
 
 A straightforward calculation yields
 \begin{multline*}
 \label{9}
 % c^{(k)}_{\lambda,n}:=
 \langle \phi_{\lambda,n}, \bar{f}_k\rangle=\int_{-\infty}^\infty
    \frac{h_k(t)e^{-2i\pi(\mu_n-\mu_k)t}}{t-(\lambda+iw)}dt \\
               =  \begin{cases}
                   e^{-2i\pi \mu_n\lambda} h_k(\lambda+iw)e^{2i\pi \mu_k \lambda}e^{2\pi w (\mu_n-\mu_k)} , & k \geq n; \\
                     0, & k< n,                  
              \end{cases}
\end{multline*}  
and
\beqq
\label{10} 
c_{\lambda,n}=e^{-2i \pi \mu_n \lambda}\sum_{k\geq n}h_k(\lambda+iw)e^{2i\pi \mu_k \lambda}e^{2\pi w(\mu_n-\mu_k)}. 
\eeqq        
 Let $d_{\lambda,n}:=c_{\lambda,n}e^{2i \pi \mu_n \lambda}$ and 
\begin{align*}
{\bf c}_{n}:=\{c_{\lambda,n}\}_{\lambda\in \Lambda}\in l^2(\Lambda), \ 
{\bf d}_{n}:=\{d_{\lambda,n}\}_{\lambda\in \Lambda}\in l^2(\Lambda), \\
          {\mathbf \omega}_\lambda=\{\omega_{\lambda, k}\}_{k \in \ZZ}:=
           \{h_k(\lambda+iw)e^{2\pi i\lambda k}\}_{k\in \ZZ}\in l^2(\ZZ).
\end{align*}
We have
\begin{align*}
%\label{11}
\|\{c_{\lambda, n}\}\|^2_{l^2(\Lambda\times \ZZ)}= \sum_n \|{\bf c}_n\|^2=\sum_n \|{\bf d}_n\|^2, \\
%\eeq
%\beq
\sum_\lambda \|\omega_\lambda\|^2=\sum_{\lambda, k}|h_k(\lambda+i)|^2 \asymp \|f\|^2,
\end{align*}
the later follows from (\ref{7}) and the fact that $\Lambda$ is a sampling set for $\PW$.

We also have 
\beq
\label{11}
{\bf d}_n = A\omega_\lambda,
\eeq 
where the matrix $A$  is defined as
\beq
\label{11a}
A=(a_{k,n})_{k,n\in \ZZ}, \quad a_{k,n}= \begin{cases}
e^{-2\pi(\mu_k-\mu_n)w}, & k \geq n; \\
0, & k< n. 
               \end{cases}.
\eeq

Denote $\gamma_n=2\pi w (\mu_{n+1}-\mu_n)$ and consider the operator $B=(b_{p,q}): l^2(\ZZ) \to l^2(\ZZ)$
with the matrix
\beqq
b_{p,q}= \begin{cases}
                 e^{-\gamma_p}, & q=p+1; \\
                 0   & \text{otherwise}
               \end{cases}.
\eeqq 
We have
\beqq 
A=I + \sum_{j\geq 1} B^j.
\eeqq
For sufficiently large $N$ such that $\mu_{n+N}-\mu_n\geq 1$ (i.e. $\gamma_n+ \ \ldots \ +\gamma_{n+N}\geq 1$) for all $n$ we have $\|B^N\|\leq e^{-\Re w}$. Therefore, this series converges to $(I-B)^{-1}$, in particular $A$ is invertible. 
%%%%%%%%%%%%%%%%%%%%%%%%%%
 %This is a Toeplitz matrix, its  symbol $(1-e^{-\pi\beta w+i\theta})^{-1}$ does not vanish for $\theta \in [-\pi,\pi]$, ,
 %$A$ is invertible. (To make the presentation self-contained, we remark that after the Fourier transform $A$ becomes multiplication by 
%$(1-e^{-\pi \beta w+i\theta})^{-1}$ in $L^2(-\pi, \pi)$).

Therefore, $\|{\bf d}_n\|\asymp \|\omega_\lambda\|$ and, finally,
\beqq
\sum_n \|{\bf c_n}\|^2 \asymp \sum_\lambda \|\omega_\lambda\|^2 \asymp \|f\|^2.
\eeqq 
This completes the proof.
% \end{proof} 

%{\em Remark.} {

{\em The statement 2) $\Rightarrow $1)} follows easily from the frame inequality. We restrict ourselves just to brief outline of its proof.
%By applying the frame inequality one can easily see that the condition for $\Lambda$ to be a sampling set in 
%$PW_{\beta(M)}$ is also necessary for $\Ga$ to be a frame in $\Lt$.  We formulate the corresponding statement 
%as a separate proposition and restrict ourselves just to brief outline of its proof.

%}

%{\bf Proposition}  
%{\em Let the sets  $\Lambda, \ M\subset \RR$    be such that the system $\Ga$ is a frame in $\Lt$. 
%Then $M$ is locally finite and also $\Lambda$ is a sampling set for $\PWm$.}
 
 %\begin{proof} 
 Indeed, relation \eqref{4b} follows from the right-hand side of \eqref{4}. In case \eqref{4b} fails
the lattice $\Lambda\times M$ contains arbitrary large clusters of points. Taking a function whose time-frequency support  is located  (mainly) in such cluster we come to contradiction with \eqref{4}.

Relation \eqref{4a} is  straightforward. If the set $M$ contains arbitrary large gaps, then  
the time-frequency shifts of the Gaussian,  centred to  the middle of these gaps bring contradiction to \eqref{4}. 

As in the proof of the previous  part  it now follows  that the operator $A$ defined in \eqref{11a}, is invertible.

In order to verify that $\Lambda$ is a sampling set in $\PWm$ it suffices to check that $\Lambda$ is a sampling set in each $PW_{[0, \beta_n]}$
with the sampling constants independent of $n$. The latter is also straightforward. Indeed, take $n\in \ZZ$, $h\in PW_{\beta_n}$, and  set
$f(t)=h(t)e^{2i\pi\mu_n t}$  (we keep the same notation as in the first part of the Theorem).  The sampling property  follows now from 
representation \eqref{11}, invertibility of the operator $A$, and also from the frame inequality. 
%\end{proof}

  \medskip
  
  {\bf Remark} {\em It  is worth to mention that, in contrast to the case of the Gaussian window, the condition for $\Ga$ to be a frame is essentially
  non-symmetric with respect to $\Lambda$ and $M$. This reflects the asymmetry of the Cauchy kernel under the Fourier transform.}

 \end{document}